\documentstyle[amssymb,amsfonts,12pt]{amsart}

\newtheorem{theorem}{Theorem}[section]
\newtheorem{proposition}[theorem]{Proposition}

\newtheorem{corollary}[theorem]{Corollary}
\newtheorem{definition}[theorem]{Definition}

\def\C{{\mbox{\rm\kern.24em
\vrule width.03em height1.43ex depth-.052ex \kern-.26em C}}}
\def\Z{{\bf Z}}
\def\R{{\mbox{\rm I\kern-.22em R}}}

\def\BMO{{\rm BMO}}

\def\X{K}  

\def\fv{{\vec f}}
\def\kv{{\vec j}}  

\def\111{\gamma}

\def\be#1{\begin{equation}\label{#1}}
\def\bas{\begin{align*}}
\def\eas{\end{align*}}
\def\bi{\begin{itemize}}
\def\ei{\end{itemize}}
\newenvironment{proof}{\noindent {\bf Proof} }{\endprf\par}
\def \endprf{\hfill  {\vrule height6pt width6pt depth0pt}\medskip}
\def\emph#1{{\it #1}}

\title{Uniform estimates on paraproducts}

\author{Camil Muscalu}
\address{Department of Mathematics, UCLA, Los Angeles CA 90095-1555}
\email{camil@@math.ucla.edu}

\author{Terence Tao}
\address{Department of Mathematics, UCLA, Los Angeles CA 90095-1555}
\email{tao@@math.ucla.edu}

\author{Christoph Thiele}
\address{Department of Mathematics, UCLA, Los Angeles CA 90095-1555}
\email{thiele@@math.ucla.edu}

\dedicatory{Dedicated to our late friend and admired mathematician Tom Wolff}

\begin{document}

\begin{abstract}  
We prove uniform $L^p$ estimates (Theorem \ref{paraproduct})
for a family of paraproducts and corresponding maximal operators.
\end{abstract}

\maketitle

\section{Introduction}\label{intro-sec}

The purpose of this article is to prove certain
uniform estimates on paraproducts. Our motivation
for this result is that we need it in a sequel \cite{mtt-3} of this
paper, where we study uniform estimates for multilinear
singular integrals as in \cite{mtt-1} with a one parameter modulation 
symmetry. Despite our current narrow objective, the uniform paraproduct 
estimates in the current article may be of independent interest 
by themselves.

As motivation, we begin by defining a 
standard paraproduct form of degree $n\ge 2$ to be 
a multilinear form of the type
\begin{equation}\label{defpara}
\Lambda(f_1,\dots, f_n)=
\int_\R \sum_{k\in \Z} \prod_{i=1}^n Q_{i,k} f_i (x)\, dx\ ,
\end{equation}
where the $Q_{i,k}$ are convolution operators
\[Q_{i,k}f=f*\phi_{i,k}\]
with convolution kernels whose Fourier transforms
$\widehat{\phi}_{i,k}$ are bump functions adapted
to $\{\xi:|\xi|\le 2^k\}$ and we assume
\begin{equation}\label{twozero}
\widehat{\phi}_{1,k}(0)=\widehat{\phi}_{2,k}(0)=0
\end{equation}
for all $k$. Here and in the sequel, a bump function $\phi$ adapted to an 
interval $I$ is a function supported in $I$ and satisfying
\[\left\|\phi^{(\alpha)}\right\|_\infty\le |I|^{-\alpha}\]
for all derivatives $\phi^{(\alpha)}$ of order 
$0\le \alpha\le N$ for some large $N$ (see \cite{stein}).

Classical Calderon-Zygmund theory gives the following standard estimate:
\begin{equation}\label{standard}
\Lambda(f_1,\dots, f_n)\le C_{p}
\prod_{i=1}^n \|f_i\|_{p_i}\end{equation}
where $1< p_i< \infty$ are any exponents satisfying the homogeneity condition
\begin{equation}\label{scaling}
\sum_{i=1}^n \frac 1{p_i}=1\ .
\end{equation}
and $C_{p}$ is a constant depending only on the tuple $(p_i)$.

To see this estimate we use H\"older's inequality to estimate the left hand side of \eqref{standard} by:
\[
\int \prod_{i=1}^2 
\left(\sum_{k=1}^\infty |Q_{ik}f_i(x)|^2\right)^{\frac 12}
\prod_{i>2}\sup_{k\in \Z} |Q_{i,k}f_i (x)|\, dx\]
\[\le \prod_{i=1}^2 \|S f_i\|_{p_i}\prod_{i>2}\|Mf_i\|_{p_i}\ \ \ .\]
Here $S$ denotes an appropriate Littlewood-Paley square function operator,
and $M$ denotes the Hardy-Littlewood maximal operator. Now the above
claim follows easily from boundedness of $S$ and $M$ in $L^p$ for $1<p<\infty$.
Observe that the same proof works if we replace
$Q_{i,k}$ by $Q_{i,k+k_i}$ in \eqref{defpara} for some
arbitrary integers $k_i$.

It is also well known (\cite{coifmanm6}) that it is sufficient to 
require condition \eqref{twozero}  only for $\widehat{\phi}_{1,k}$.
We shall sketch a proof under the more restrictive assumption
that $\widehat{\phi}_{1,k}$ vanishes of $\{\xi:|\xi|\le 2^{k-L}\}$
for some given $L\ge 1$.
Indeed, in this case we can write for each $i>1$
\[Q_{i,k}=\sum_{l=1}^{n+L} Q_{i,k,l}\ \ \ ,\]
where the symbol $\widehat{\phi}_{i,k,l}$
of the Fourier multiplier $Q_{i,k,l}$, is adapted to 
$\{\xi:|\xi|\le 2^{k-l}\}$ and vanishes on the
smaller interval 
$\{\xi:|\xi|\le 2^{k-l-2}\}$ except in the case $l=n+L$ when
we do not require any vanishing.
Then we can split $\Lambda(f_1,\dots, f_n)$ as in \eqref{defpara}
by the distributive law
into a sum of polynomially in $L$ many multilinear operators, all of which are multiples of standard
paraproducts after possibly permuting the indices, 
except for the term
\[\int \sum_k Q_{1,k}f_1(x)\prod_{i>1} Q_{i,k,n+L}f_i(x) \, dx\ .\]
However, this term simply vanishes, because the Fourier transform
of
\[\prod_{i>1} (f_i*\phi_{i,k,n+L})\]
is supported in $\{\xi:|\xi|\le  n2^{k-n-L}\}$ whereas the Fourier transform
of 
\[\overline{f_1*\phi_{1,k}}\]
vanishes on $\{\xi:|\xi|\le 2^{k-L}\}$, so the inner product of these two vanishes
by Plancherel.

Observe that in the above argument, the constant in the estimate for
$\Lambda$ depends on $L$. However, this dependence is artificial
as one can show a bound \eqref{standard} with a constant 
independent of $L$ and independent of any constants $k_i$ as above
that we may insert in the definition of $\Lambda$.
The purpose of this article is to prove a very 
general version of this uniformity result:

\begin{theorem}\label{paraproduct}  Let $n \geq 2$, and let $\Omega \in \Z^n$ be a set of the form
\begin{equation}\label{omegaconst}
 \Omega = \{ (j_1, \ldots, j_n): j_{i_\alpha} \geq j_{i'_\alpha} + A_\alpha \hbox{ for all } 1\le \alpha \le \X\}
\end{equation}
where $\X\ge 0$ is an integer and $i_\alpha, i'_\alpha \in \{1, \ldots, n\}$, $A_\alpha \in \Z$ for all $1 \leq \alpha \leq \X$.  For each $1 \leq i \leq n$ and $j \in \Z$, let $Q_{j,i}$ be a Fourier multiplier whose symbol is a bump function adapted to $\{ \xi: |\xi| \le 2^j \}$ that vanishes at the origin.  Then one has the estimates
\be{para-est}
\left|\int_\R \sum_{(j_1, \ldots, j_n) \in \Omega} \prod_{i=1}^n Q_{j_i,i} f_i\right|
\leq C_{p,n,K} \prod_{i=1}^n \|f_i\|_{p_i}
\end{equation}

\be{max-para-est}
\left\|\sup_{j_n\in \Z} \left|\sum_{(j_1,\dots,j_n)\in \Omega} \prod_{i=1}^{n-1}Q_{j_i,i}f_i\right|
\right\|_{p_n'}
\le C_{p,n,K} \prod_{i=1}^{n-1}\|f_i\|_{p_i}
\end{equation}
for all $1 < p_i < \infty$ obeying 
\eqref{scaling}, where $1/p_n'+1/p_n=1$
and the constant $C_{p,n,K}$ depends on the $p_i$, $n$, $K$, 
but is independent of the $A_i$ and $f_i$.
\end{theorem}

The set $\Omega$ in this theorem is a convex
polytope in $\Z^n$ which is invariant under translation in direction 
$(1,\dots,1)$. It has essentially $K$ faces.
We remark that the number $K$ of independent constraints in 
\eqref{omegaconst} can be bounded by a function of $n$, so the dependence of
$C_{p,n,K}$ on $K$ serves only to emphasize that we have 
better estimates for small $K$.

To compare this result further with the existing literature, we observe that
we can write a paraproduct as in \eqref{defpara} or \eqref{para-est}
in the multiplier form
$$\int_{\sum\xi_i=0} m(\xi_1,\dots,\xi_n) \prod_{i=1}^n \widehat{f}_i(\xi_i)\,d \sigma$$
where the multiplier $m$ satisfies the symbol estimates
\begin{equation}\label{symbest}
\partial^\alpha m(\xi)\le C |\xi|^{-|\alpha|}
\end{equation}
for any multi index $\alpha$ up to some fixed large order
and for all vectors $\xi$. Under this condition on the multiplier,
standard paraproduct theory (e.g. \cite{coifmanm1}-\cite{coifmanm6}) gives 
\eqref{para-est}.
However, for the multipliers arising from Theorem \ref{paraproduct}, 
the constant in \eqref{symbest}
depends on the set $\Omega$. Hence the point of
the theorem is uniformity in $\Omega$.

The multipliers of the theorem do satisfy a bound of the type
\begin{equation}\label{marc}
(\prod_{i=1}^n \partial_i^{\alpha_i}) m(\xi)
\le C\prod_{i=1}^n |\xi_i|^{-\alpha_i}\ \ \ .
\end{equation}
with a constant $C$ independent of $\Omega$.
However, as is shown in \cite{graf-kalt}, this condition is in general
not sufficient to guarantee an estimate \eqref{para-est}.
If the multiplier $m$ is of tensor product type
\[m(\xi)=\prod_{i=1}^n m_i(\xi_i)\]
in addition to satisfying \eqref{marc}, then the integrand in the definition
of $\Lambda$ splits into a tensor product and one can trivially show
\eqref{para-est} by reducing to the one dimensional case.
Observe that the multipliers arising in Theorem \ref{paraproduct}
do not split as tensor products because of the constraints on
$\Omega$. However, the idea of the theorem is that the special structure
of $\Omega$ still allows to recover enough of the good behaviour of the 
tensor product case.

We shall not persue here results for endpoints 
$p_i = \infty$ or $p_i = \BMO$.

We prove Theorem \ref{paraproduct} in Sections \ref{notation1-sec}-\ref{dominant-sec}. It is obtained by induction on the dimension $n$ of $\Omega$.  
The main tools are the H\"older, Littlewood-Paley, and Hardy-Littlewood maximal inequalities, together with some combinatorial manipulations to break $\Omega$ up into simpler objects.  The latter shall be most easily accomplished using the language of graph theory. 

We shall specialize Theorem \ref{paraproduct} to a form which will be convenient for our application to multilinear singular integrals with modulation
symmetries.

\begin{corollary}\label{paraproduct-cor}  Let $n \geq 2$, and let $M_1, \ldots, M_n$ be integers.  For each $1 \leq i \leq n$ and $j \in \Z$, let $\pi_{j,i}$ be a Fourier multiplier whose symbol is a bump function adapted to $\{ \xi: |\xi| \le 2^{j + M_i} \}$.  Suppose that for each $k \in \Z$ there exists at least one $1 \leq i \leq n$ such that the symbol of $\pi_{j,i}$ vanishes at the origin. Then one has the estimate
\be{para-cor}
\sum_j |\int \prod_{i=1}^n \pi_{j,i} f_i|
\leq C_{p,n,K} \prod_{i=1}^n \|f_i\|_{p_i}
\end{equation}
for all $1 < p_i < \infty$ obeying 
$\sum_{i=1}^n \frac 1{p_i}=1$, 
where the constant $C_{p,n,K}$ depends on the $p_i$, $n$, $K$,
but is independent of the $M_i$ and $f_i$.
\end{corollary}

We prove this Corollary in Section \ref{cor-proof}.  It shall follow easily from Theorem \ref{paraproduct} after some normalizations and dyadic decompositions on the $\pi_{j,i}$.  It is easy to see that some sort of cancellation condition on the $\pi_{j,i}$, such as the vanishing condition given above, is necessary, as can be seen by testing \eqref{para-cor} when the $\hat f_i$ are all equal to the same approximation to the delta function at the origin.  It is important that the absolute values in \eqref{para-cor} are not inside the integral, as the estimate is easily seen to be false otherwise.

The first author was partially supported by a Sloan Dissertation Fellowship.
The second author is a Clay Prize Fellow and is supported by grants from
the Sloan and Packard Foundations. The third author was partially supported
by a Sloan Fellowship and by NSF grants DMS 9985572 and DMS 9970469.

We would like to point out that Tom Wolff's first major
result, his simplification of the proof of the corona theorem,
heavily relied on estimates for paraproducts.

\section{Notation and preliminaries}\label{notation1-sec}

It will be convenient to write the paraproducts using the language of graph theory.

\begin{definition}\label{graph-def}  Let $V$ be a finite set.  A \emph{(weighted directed) edge} supported in $V$ is a triplet $e = (v_1, v_2, A) \in V \times V \times \Z$.  We refer to $v_1 =: v_1(e)$, $v_2 =: v_2(e)$, and $A(e)$ as the \emph{initial vertex}, \emph{final vertex}, and \emph{weight} of $e$ respectively.
A \emph{(weighted directed) graph} $G = (V, E)$ is a finite collection of vertices $V$ and  a finite collection $E$ of weighted directed edges $e$ supported in $V$.  We call $|V|$ and $|E|$ the \emph{order} and \emph{complexity} of $G$ respectively.
\end{definition}

We allow multiple edges from one vertex to another, as well as edges from a vertex to itself, although such edges are usually redundant or trivial in our applications.  We shall henceforth omit the modifiers ``weighted directed'' in the sequel.

If $p$ is a directed path in $G$, we define the \emph{weight} $A(p)$ of $p$ to be the sum of all the weights of the edges in $G$.  If $V'$ is a subset of  $G$, we define $G|_{V'}$ to be the graph with vertex set $V'$ and edges $\{ e \in E: v_1(e), v_2(e) \in V'\}$.  We write $G - V'$ for $G|_{V \backslash V'}$.

If $E'$ is a set of edges whose vertices are in $V$, we write $G \cup E'$ for $(V, E \cup E')$ and $G - E'$ for $(V, E - E')$.

If $G = (V,E)$ is a graph and $V_1$, $V_2$ are disjoint subsets of $G$, we say that $V_1$, $V_2$ are \emph{adjacent} if there exists an edge in $G$ with one vertex in $V_1$ and the other in $V_2$.

Every graph can be associated with a polytope:

\begin{definition}\label{polytope-def}  Let $G = (V, E)$.  We define the \emph{polytope} $\Omega(G) \subset \Z^V$ of $G$ to be the set
$$ \Omega(G) := \{ (j_v)_{v \in V} \in \Z^V: j_{v_1(e)} \geq j_{v_2(e)} + A(e) \hbox{ for all } e \in E \}.$$
We shall often write $\kv$ for $(j_v)_{v \in V}$.
\end{definition}

Note that $\Omega(G)$ is empty if $G$ contains directed circuits of positive weight.  (This implication can be reversed, but we shall not need this here).  Conversely, if $E$ is empty, then $\Omega(G) = \Z^V$.

\begin{definition}  If $e = (v_1, v_2, A)$ is an edge, we define the \emph{dual} $e^*$ of $e$ by $e^* := (v_2, v_1, -A+1)$.
\end{definition}

We observe the basic relationship
\be{split}
\Omega(V, E) = \Omega(V, E \cup \{e\}) \uplus \Omega(V, E \cup \{e^*\})
\end{equation}
whenever $(V,E)$ is a graph and $e$ is an edge with vertices in $V$, and $\uplus$ denotes disjoint union.  The identity \eqref{split} allows us to replace an edge $e$ in a graph by its dual, modulo objects of lesser complexity.

Theorem \ref{paraproduct} is concerned with a certain graph of order $n$. We shall fix this
$n$: in what follows all graphs shall have order less than $n$. Observe that while
we do not exclude multiple edges, in what follows we can easily reduce graphs to not
contain multiple edges, and thus the complexity of all graphs can be
assumed to be less than $(n+1)^2$.
For each $v\in V$ and $j \in \Z$, we fix $Q_j^v$ to be a Fourier multiplier whose symbol is a bump function adapted to $\{ \xi: 2^{10nj} \leq |\xi| \leq 2^{10nj+1} \}$.  We use powers of $2^{10n}$ rather than $2$ to obtain better lacunary separation properties.

We shall need two types of paraproducts: non-maximal paraproducts, which we define as multi-linear forms; and maximal paraproducts, which we define as multi-linear operators.  

\begin{definition}\label{paraproduct-def}
Let $G = (V,E)$ be a graph with $|V|\le n$, and suppose that $\fv = (f_v)_{v \in V}$ is a family of test functions on $\R$.  We define the \emph{non-maximal paraproduct} $\Lambda_G$ to be the quantity
$$ \Lambda_G( \fv ) :=
\sum_{\kv \in \Omega(G)} \int \prod_{v \in V} Q_{j_v}^v f_v.$$
We define the \emph{order} and \emph{complexity} of $\Lambda_G$ to be the order and complexity of $G$ respectively. 
\end{definition}

\begin{definition}\label{max-paraproduct-def}
Let $G = (V \uplus \{v_*\},E)$ be a graph with a distinguished vertex $v_*$ and $|V|\le n$, 
and suppose $\fv = (f_v)_{v \in V}$ is a family of test functions on $\R$.  We define the \emph{maximal paraproduct} $\Lambda^{[v_*]}_G$ to be the function
$$ \Lambda^{[v_*]}_G( \fv )
= \sup_{j_* \in \Z} |\sum_{\kv \in \Omega(G): j_{v_*} = j_*} \prod_{v \in V} Q_{j_v}^v f_v|.$$
We define the \emph{order} and \emph{complexity} of $\Lambda^{[v_*]}_G$ to be $|V| + \frac{1}{2}$ and $|E|$ respectively.
\end{definition}

Note that a maximal paraproduct on $n$ functions is slightly higher order than a non-maximal paraproduct on $n$ functions, but slightly lower order than a non-maximal paraproduct on $n+1$ functions.  From \eqref{split} we observe the identity
\be{para-split}
\Lambda_{G}(\fv) = \Lambda_{G \cup \{e\}}(\fv) + \Lambda_{G \cup \{e^*\}}(\fv)
\end{equation}
and the triangle inequalities
\be{max-split}
\Lambda^{[v_*]}_{G}(\fv)
\leq
\Lambda^{[v_*]}_{G \cup \{e\}}(\fv) + \Lambda^{[v_*]}_{G \cup \{e^*\}}(\fv)
\end{equation}
and
\be{max-split-alt}
\Lambda^{[v_*]}_{G \cup \{e\}}(\fv) \leq
\Lambda^{[v_*]}_{G}(\fv) + \Lambda^{[v_*]}_{G \cup \{e^*\}}(\fv)
\end{equation}
whenever $e$ is an edge with vertices in the vertex set of $G$.

We shall prove 

\begin{theorem}\label{paraproduct-graph}  
\begin{itemize}
\item Let $G = (V,E)$ be a non-empty graph with $|V|\le n$, and $\fv=(f_v)_{v \in V}$ a family of test functions.  Then we have
\be{lag}
|\Lambda_G( \fv )| \leq C_{(p_v),n} \prod_{v \in V} \|f_v\|_{p_v}
\end{equation}
whenever $1 < p_v < \infty$ and $\sum_{v \in V} 1/p_v = 1$.
\item Let $G = (V \uplus \{v_*\},E)$ be a non-empty 
graph with $|V|\le n$, and $\fv=(f_v)_{v \in V}$ a family of test functions.  Then we have
\be{lag-max}
\|\Lambda^{[v_*]}_G( \fv )\|_{p_*} \leq C_{p_*,(p_v),n} \prod_{v \in V} \|f_v\|_{p_v}
\end{equation}
whenever $1 < p_*, p_v < \infty$ and $\sum_{v \in V} 1/p_v = 1/p_*$.
\end{itemize}
\end{theorem}

Theorem \ref{paraproduct} then follows easily from rescaled versions of \eqref{lag} and a dyadic decomposition of the multipliers $Q_{k,j}$.

The theorem is easily verified from the standard linear theory of maximal truncated singular integrals (see e.g. \cite{stein}) when the order of the paraproduct is $1\frac{1}{2}$ (the lowest order in which the statement is not void).  
Now suppose the order is at least two and let $n=|V|$ be the number of functions. 
We shall divide into two cases.

\begin{definition}\label{freq-def}  

\begin{itemize}
\item We say that $\Lambda_G$ (resp. $\Lambda^{[v_*]}_G$) has a \emph{dominant frequency} $v_0 \in V$ if one has 
$$ j_{v_0} > j_v$$
whenever $\kv \in \Omega(G)$ and $v \in V \backslash \{v_0\}$.
\item
We say that $\Lambda_G$ (resp. $\Lambda^{[v_*]}_G$)
has \emph{two competing frequencies} $v_1, v_2 \in V$ if $v_1\neq v_2$ 
and 
$$ j_{v_1} = j_{v_2}$$
whenever $\kv \in \Omega(G)$.  A paraproduct with two competing frequencies $v_1, v_2$ is in \emph{standard form}
if $\{v_2\}$ is not adjacent to $V \backslash \{v_1, v_2\}$ (resp.
$ \{v_*\}\cup V \backslash \{v_1,v_2\}$).
\end{itemize}
\end{definition}

Note that the frequency $j_{v_*}$ plays no role in determining whether a maximal paraproduct has a dominant frequency or two competing frequencies.

From \eqref{para-split}, \eqref{max-split} and the introduction of edges of the form $(v,v',0)$ or their duals,
we observe that any paraproduct can be estimated by the sum of at most $C_n$ paraproducts of the same order (each with complexity increased by at most $C_n$), such that each such paraproduct either has a dominant frequency or two competing frequency.  Note that any paraproduct with two competing frequencies can be placed in standard form, by replacing any edge with vertex $v_2$ with the corresponding edge with vertex $v_1$, and then adding the edges $(v_1,v_2,0)$ and $(v_2,v_1,0)$.  

If a non-maximal paraproduct $\Lambda_G$ has a dominant frequency $v_0$, then the summands have frequency approximately $2^{10nj_{v_0}}$, and so $\Lambda_G$ vanishes identically.

Now suppose a non-maximal paraproduct $\Lambda_G$ has two competing frequencies $v_1$, $v_2$ and is in standard form.  Assuming first $V\setminus\{v_1,v_2\}\neq \emptyset$ 
we can use Cauchy-Schwarz to estimate
$$
|\sum_{\kv \in \Omega(G)} \prod_{v \in V} Q_{j_v}^v f_v|
\leq S_{v_1}(f_{v_1}) S_{v_2}(f_{v_2}) \sup_{j_*} |\sum_{\kv \in \Omega(G): j_{v_1} = j_{v_2} = j_*} \prod_{v \in V \backslash \{v_1,v_2\}} Q_{j_v}^v f_v|$$
where $S_v(f_v)$ is the Littlewood-Paley square function 
$$ S_v(f_v) := (\sum_k |Q_k^v f_v|^2)^{1/2}.$$
From this and H\"older's inequality, and the boundedness of $S_v$ on $L^p$, $1 < p < \infty$ we have
\be{compete-est}
\Lambda_G((f_v)_{v \in V}) \leq C_{p_1,p_2}
\| f_{v_1}\|_{p_1} 
\| f_{v_2} \|_{p_2} 
\| \Lambda^{[v_1]}_{G|_{V \backslash \{v_2\}}}
(\{f_v\}_{v \in V \backslash \{v_1,v_2\}}) \|_{p_*}
\end{equation}
whenever $1 < p_1, p_2, p_* < \infty$ are such that $1/p_1 + 1/p_2 + 1/p_* = 1$.
If $V\setminus\{v_1,v_2\}=\emptyset$, then the above calculation holds without
the factor involving the maximal paraproduct, and we immediately conclude
Theorem \ref{paraproduct-graph}.

From the above remarks we see that the estimate \eqref{lag} for non-maximal paraproducts 
automatically follows from the estimate \eqref{lag-max} for maximal paraproducts of smaller
order. Thus it remains to consider maximal paraproducts.

To estimate maximal paraproducts we require some further definitions.  

\begin{definition}\label{semi-direct-def}
\begin{itemize}
\item A maximal paraproduct $\Lambda^{[v_*]}_G$ is \emph{separable} if one can write $V = V_1 \uplus V_2$ where $V_1$ and $V_2$ are non-empty and not adjacent.

\item A maximal paraproduct $\Lambda^{[v_*]}_G$ with dominant frequency $v_0 \in V$ is said to be \emph{semi-direct} if $\{v_*\}$ and $V \backslash \{v_0\}$ are not adjacent.

\item A maximal paraproduct $\Lambda^{[v_*]}_G$ with two competing frequencies $v_1, v_2 \in V$ is said to be \emph{semi-direct} if $\{v_*\}$ and $V \backslash \{v_1,v_2\}$ are not adjacent.

\item A maximal paraproduct is said to be \emph{good} if it is of one of the above three types.
\end{itemize}
\end{definition}

We now show that the estimate \eqref{lag-max} for good maximal paraproducts follows from an 
application of Theorem \ref{paraproduct-graph} applied to paraproducts of strictly lower order.

If $\Lambda^{[v_*]}_G$ is separable, then we have the pointwise estimate
\be{sep-max-est}
\Lambda^{[v_*]}_G( (f_v)_{v \in V} ) \leq
\Lambda^{[v_*]}_{G|_{V_1 \cup \{v_*\}}}( (f_v)_{v \in V_1} )
\Lambda^{[v_*]}_{G|_{V_2 \cup \{v_*\}}}( (f_v)_{v \in V_2} ).
\end{equation}
The claim then follows from H\"older.

Now suppose $\Lambda^{[v_*]}_G$ is semi-direct with dominant frequency $v_0$. We shall
assume first $V\setminus \{v_0\}$ is not empty.
 In this case we observe the identity
$$ 
\sum_{\kv \in \Omega(G): j_{v_*} = j_*} \prod_{v \in V} Q_{j_v}^v f_v
= \sum_{j_0:(j_*,j_0) \in \Omega(G|_{\{v_*,v_0\}})}\ 
\sum_{\kv \in \Omega(G - \{v_*\}): j_{v_0} = j_0} \prod_{v \in V} Q_{j_v}^v f_v
$$
for all $j_*$.  We may of course assume that $\Omega(G|_{\{v_*,v_0\}})$ is non-empty.  From the dominant frequency hypothesis we thus have $j_{v_0} > j_v$ whenever
$\kv \in \Omega(G - \{v_*\})$ and $v \in V \backslash \{v_0\}$.  In particular, the inner summand has Fourier transform supported on the annulus $\{ 2^{10nj_0-3n} \leq |\xi| \leq 2^{10nj_0+3n} \}$.  Since 
$\{j_0:(j_*,j_0)\in \Omega(G|_{\{v_*,v_0\}})\}$ is an interval (possibly infinite)
for each $v_*$, we may omit the constraint on $j_0$ and replace it by an appropriate Fourier 
multiplier cutting of the very high and very low frequencies (not needed if the interval
is infinite on either side). This Fourier multiplier can be estimated
by the Hardy Littlewood maximal function. Thus we have the pointwise estimate
$$
|\sum_{\kv \in \Omega(G): j_{v_*} = j_*} \prod_{v \in V} Q_{j_v}^v f_v|
\leq C M |\sum_{j_0 \in \Z}
\sum_{\kv \in \Omega(G - \{v_*\}): j_{v_0} = j_0} \prod_{v \in V} Q_{j_v}^v f_v|$$
where $M$ is the Hardy-Littlewood maximal operator.  Taking suprema in $j_*$ and applying the Hardy-Littlewood maximal and Littlewood-Paley  inequalities,
we obtain
$$ 
\| \Lambda^{[v_*]}_G( (f_v)_{v \in V} )\|_{p_*} \leq
C \| (\sum_{j_0 \in \Z} |\sum_{\kv \in \Omega(G - \{v_*\}): j_{v_0} = j_0} \prod_{v \in V} Q_{j_v}^v f_v|^2)^{1/2} \|_{p_*}.$$
The expression inside the norm on the right-hand side is pointwise bounded by
$$ S_{v_0}(f_{v_0}) \Lambda^{[v_0]}_{G - \{v_*\}}( 
(f_v)_{v \in V \backslash \{v_0\}} ).$$
The claim then follows by H\"older.
The case $V\setminus\{v_0\}=\emptyset$ is again an easy variant of the
above without a smaller order maximal paraproduct.

Finally, suppose $\Lambda^{[v_*]}_G$ is semi-direct with two competing frequencies $v_1, v_2 \in V$.  We may place this paraproduct in standard form without affecting the property of being semi-direct.  But we may then repeat the arguments used to prove \eqref{compete-est}, and obtain the pointwise estimate
$$
\Lambda^{[v_*]}_G((f_v)_{v \in V}) \leq 
S_{v_1}(f_{v_1}) S_{v_2}(f_{v_2}) \Lambda^{[v_1]}_{G - \{v_*,v_2\}}(\{f_v\}_{v \in V \backslash \{v_1,v_2\}}),
$$
again with no maximal paraproduct on the right hand side if $V\setminus\{v_1,v_2\}$ is empty.
The claim then follows from H\"older and the boundedness of the Littlewood-Paley square function on $L^p$, $1 < p < \infty$.

In light of the preceding discussion and induction, Theorem \ref{paraproduct-graph} (and hence Theorem \ref{paraproduct}) will follow if we can show

\begin{proposition}\label{decomp}  Every maximal paraproduct $\Lambda^{[v_*]}_G$ can be bounded by a finite collection of good maximal paraproducts of the same order.  The cardinality of the collection is bounded by a quantity depending only on the order and complexity of $\Lambda^{[v_*]}_G$.
\end{proposition}

This will be done in the next two sections.  The basic idea shall be to exploit variants of the identity
$$ \chi_{a < b, a< c} = \chi_{a < b \leq c} + \chi_{a < c < b}$$
as well as \eqref{max-split}, \eqref{max-split-alt} to permute the edges of the graph $G$ into a good form. 

\section{The case of two competing frequencies}\label{compete-sec}

It remains only to prove Proposition \ref{decomp}.  By the remarks in the previous section it suffices to verify the Proposition when $\Lambda^{[v_*]}_G$ has a dominant frequency or two competing frequencies.

In this section we consider the case of two competing frequencies $v_1$, $v_2$; this case is less technical than the dominant frequency case.  We may place the maximal paraproduct in standard form.  The vertex $v_2$ now plays no significant role, and we shall often work with the reduced graph $G - \{v_2\}$.  Let $E - \{v_2\}$ denote the edge set of $G - \{v_2\}$.

Write $G = (V \cup \{v_*\}, E)$.  We shall induct on the complexity of $G$.  Suppose inductively that the claim is proven for all maximal paraproducts $\Lambda^{[v_*]}_{(V \cup \{v_*\},E')}$
with two competing frequencies $v_1, v_2$ in normal form, with $|E'| < |E|$.

In light of \eqref{max-split-alt} we may replace any edge $e$ in $E - \{v_2\}$ by its dual $e^*$ modulo objects of lower complexity, without disturbing the property of having two competing frequencies in standard form.  Thus we may freely reverse the direction of edges in $G - \{v_2\}$ (adjusting the weights accordingly).

We may assume that $G - \{v_2\}$ is connected as an undirected graph, for if $G - \{v_2\}$ split into two disconnected components then $\Lambda^{[v_*]}_G$ is clearly good.

We may also assume that $G - \{v_2\}$ contains no circuits as an undirected graph.  To see this, we suppose for contradiction that $G - \{v_2\}$ contained an undirected circuit.  By the freedom to reverse the direction of edges in $G - \{v_2\}$ we may assume that this circuit is directed, and has positive weight.  But then $\Lambda^{[v_*]}_G$ is identically zero, and we are done.

From the preceding discussion we see that we may assume that $G - \{v_2\}$ is a tree (when viewed as an undirected graph), with no loops or double edges.  We may then assume that the vertex $v_*$ is a leaf node (i.e. has degree exactly one), since $\Lambda^{[v_*]}_G$ becomes separable otherwise.

We now induct on the distance from $v_*$ to $v_1$ in $G - \{v_2\}$ (viewed as an undirected graph).  If this distance is 1 then $\Lambda^{[v_*]}_G$ is semi-direct and we are done.  Now suppose inductively that the distance is greater than 1.  Then $v_*$ is adjacent to some vertex $v' \in V - \{v_1,v_2\}$, which is in turn adjacent to some vertex $v'' \in V - \{v_2\}$, such that $v''$ is closer to $v_1$ than $v'$ in $G - \{v_2\}$ (viewed as an undirected graph).

By reversing edges if necessary we may assume that $G$ contains the edges $e_1 = (v_*, v', A_1)$ and $e_2 = (v', v'', A_2)$.  Observe that the statement
$$ j_{v_*} \geq j_{v'} + A_1 \hbox{ and } j_{v_*} \geq j_{v''} + A_1 + A_2$$
holds if and only if one of the (exclusive) statements
$$ j_{v_*} \geq j_{v'} + A_1 \hbox{ and } j_{v'} \geq j_{v''} + A_2$$
or
$$ j_{v_*} \geq j_{v''} + A_1 + A_2 \hbox{ and } j_{v''} \geq j_{v'} - A_2 + 1$$
hold.  From this and the triangle inequality we have
\be{sep-chain}
\Lambda^{[v_*]}_G(\fv) \leq \Lambda^{[v_*]}_{G - \{e_2\} \cup \{e_3\}}(\fv)
+ \Lambda^{[v_*]}_{G - \{e_1, e_2\} \cup \{e_2^*, e_3\}}(\fv)
\end{equation}
where $e_3 := (v_*, v'', A_1 + A_2)$.  The first term on the right-hand side is separable, while the second term on the right-hand side can be treated by the induction hypothesis.  This completes the inductive argument.

\section{The case of a dominant frequency}\label{dominant-sec}

To conclude the proof of Proposition \ref{decomp} and hence Theorem \ref{paraproduct} it remains only to consider the case when $G$ has a dominant frequency $v_0$.  In this case we may freely replace $G$ by $G \cup H$, where $H$ is the graph with vertex set $V \cup \{v_*\}$ and edge set
$$ \{ (v_0, v, 1): v \in V - \{v_0\} \},$$
since this makes no difference to the maximal paraproduct $\Lambda^{[v_*]}_G$ other than to increase the complexity slightly.

Fix $V$, $v_0$, $v_*$.  The graph $H$ is an example of a \emph{heirachy}, which we now define.

\begin{definition}\label{heirarchy-def}  A \emph{heirarchy} is a graph $H = (V \cup \{v_*\}, E')$ obeying the following properties.
\begin{itemize}
\item No edge of $E'$ has $v_*$ as a vertex.
\item For every $v \in V - \{v_0\}$ there exists a unique directed path $p(v) = (e_1, e_2, \ldots, e_r)$ of edges in $E'$ from $v_0$ to $v$.   We refer to $r$, $e_r$, and $A(p(v))$ as the \emph{ranking}, \emph{link}, and \emph{depth} of $v$ respectively, and define the \emph{superior} $s(v)$ of $v$ to be the initial vertex of $e_r$.  We define $v_0$ to have 0 ranking and depth, and no link or superior.
\item We have $A(p(v)) > 0$ for all $v \in V - \{v_0\}$.  
\end{itemize}
\end{definition}

Note that if $H$ is a heirarchy and $G$ is any graph on $V \cup \{v_*\}$, then $\Lambda^{[v_*]}_{G \cup H}$ automatically has $v_0$ as a dominant frequency.
Indeed one has 
\be{chunk}
j_v \leq j_{v_0} - A(p(v)) 
\end{equation}
for all $v \in V$ and $\kv \in \Omega(H)$.  Also observe that $H - \{v_*\}$ is necessarily a rooted tree with root $v_0$.

Proposition \ref{decomp} now clearly follows from

\begin{proposition}\label{decomp-dominant}  Let $G$ be a graph on $V \cup \{v_*\}$ and let $H$ be a heirarchy.  Then $\Lambda^{[v_*]}_{G \cup H}$ can be bounded by finitely many good maximal paraproducts.  The number of such paraproducts depends only on the complexity and order of $G$.
\end{proposition}

\begin{proof}
As before, we induct on the complexity of $G$.  We may assume that the proposition is true for all graphs $G$ of lower complexity and all heirarchies $H$.
The inductive hypothesis is allowed to be vacuous when the complexity 
of $G$ is zero.

By the arguments in the preceding section, we may freely reverse edges in $G$ (but not in $H$!), and assume that $G$ contains no circuits as an undirected graph.

We may assume that $G$ does not contain any edge which has the same initial and final vertices as one in $H$, since the edge from $G$ is either redundant or can be used to replace the edge in $H$, and in either case we can reduce the complexity of $G$.

We divide into two cases, depending on whether $G - \{v_*\}$ is non-empty or not.

{\bf Case 1.}  Suppose that $G - \{v_*\}$ is non-empty.  Then $G$ contains an edge $e = (v,v',A_0)$ for some $v, v' \in V$.  By replacing $e$ with its dual if necessary we may assume that 
\be{chain}
A_0+A(p(v)) > A(p(v'))\  .
\end{equation}

If $v' \in p(v)$ (i.e., the path $p(v)$ crosses the vertex $v'$,) then the 
condition \eqref{chain} and the fact that $H$ is a heirarchy forces $G \cup H$ to contain a directed circuit of positive weight, so that $\Lambda^{[v_*]}_{G \cup H}$ vanishes.  Hence we may assume that $v' \not \in p(v)$.  In particular we have $v' \neq v_0$.

Let $e_1$ be the link of $v'$.  From \eqref{max-split-alt} we have
$$ \Lambda^{[v_*]}_{G \cup H}(\fv) \leq
\Lambda^{[v_*]}_{G \cup H - \{e_1\}}(\fv)
+ \Lambda^{[v_*]}_{G \cup H - \{e_1\} \cup \{e_1^*\}}(\fv).
$$
The graph $G \cup H - \{e_1\}$ can be written as the union of $G - \{e\}$
and $H - \{e_1\} \cup \{e\}$.  Since $H$ is a heirarchy, one can verify using \eqref{chain} that $H - \{e_1\} \cup \{e\}$ is also a heirarchy (in fact, the depth of any vertex either stays constant or increases).  Thus this summand is acceptable by the induction hypothesis.

Now consider the contribution of $G \cup H - \{e_1\} \cup \{e_1^*\}$.  If the superior $s(v')$ of $v'$ is equal to $v_0$ then one can see using \eqref{chain} that this graph contains a directed circuit of positive weight, so the contribution of this graph vanishes.  Now suppose $s(v') \neq v_0$.  Let $e_2$ be the link of $s(v')$. 

Observe from \eqref{chunk}, \eqref{chain} that the edge $e_2$ is redundant in $G \cup H - \{e_1\} \cup \{e_1^*\}$:
$$ \Omega(G \cup H - \{e_1\} \cup \{e_1^*\}) = \Omega(G \cup H - \{e_1, e_2\} \cup \{e_1^*\}).$$
The graph $G \cup H - \{e_1, e_2\} \cup \{e_1^*\}$ can be written as the union of $G - \{e\}$ and $H - \{e_1, e_2\} \cup \{e, e_1^*\}$.  The latter graph can be verified to be a heirarchy using \eqref{chain} by similar arguments as before, and so this case is also treatable by the induction hypothesis.  This concludes the treatment of Case 1.

{\bf Case 2.}  Suppose that $G - \{v_*\}$ is empty.  Define the \emph{total ranking} to be the sum of rankings of all the vertices in $V$ which are adjacent to $v^*$.

We induct on the total ranking.  If the total ranking is zero then  $\Lambda^{[v_*]}_{G \cup H}$ is semi-direct with dominant frequency $v_0$.  Now suppose the total ranking is positive.  Then $v_*$ is adjacent to at least one vertex $v$ in $V - \{v_0\}$.  By reversing edges if necessary, we may assume that $G$ contains an edge $e$ of the form $e = (v_*, v, A_0)$.

Let $e_1 = (s(v), v, A_1)$ be the link of $v$, and let $e_2$ denote the edge 
$e_2 := (v_*, s(v), A_0 - A_1)$.  By repeating the derivation of \eqref{sep-chain} we have the triangle inequality
$$
\Lambda^{[v_*]}_{G \cup H}(\fv) \leq
\Lambda^{[v_*]}_{G - \{e\} \cup \{e_2\} \cup H}(\fv)
+
\Lambda^{[v_*]}_{G \cup H - \{e_1\} \cup \{e_2^*\}}(\fv).$$
The first paraproduct on the right-hand side is of the same form as $G \cup H$, but has the total ranking reduced by one, and so can be handled by the induction hypothesis.  Since $H$ is a tree, the second paraproduct is separable.  This concludes the treatment of Case 2.
\end{proof}

Proposition \ref{decomp} has now been proven in all cases.  The proof of Theorem \ref{paraproduct} is now complete.

\section{Proof of Corollary \ref{paraproduct-cor}}\label{cor-proof}

Fix $n$, $M_1, \ldots, M_n$, $\pi_{j,i}$ as in Corollary \ref{paraproduct-cor}.
Let $m_{j,i}$ denote the symbol of $\pi_{j,i}$.  By the hypotheses of the Corollary, pigeonholing,
and symmetry  we may assume that $m_{j,n}(0) = 0$ for all $j \in \Z$.
By shifting the summation index we may assume $M_n=0$.

We may replace the left-hand side of \eqref{para-cor} by
$$ |\sum_j \epsilon_j \int \prod_{i=1}^n \pi_{j,i} f_i|$$
where $\epsilon_j$ are bounded constants.  By absorbing the $\epsilon_j$ into $m_{j,n}$ (for instance) we may assume that $\epsilon_j = 1$.

Following \cite{thiele}, we now renormalize the $m_{j,i}$ in a form which allows a good decomposition as a telescopic series.  For each integer $j$, let $\phi_j$ be a bump function adapted to $[-2^{j}, 2^{j}]$ which satisfies $\phi_j(0)=1$.

For any $1 \leq i \leq n$, we decompose 
$$ m_{j,i} = m_{j,i}(0) \phi_{j+M_i} + (m_{j,i} - m_{j,i}(0) \phi_{j+M_i}).$$
Thus by decomposition, we may assume that
$$
m_{j,i} \hbox{ is a bounded multiple of } \phi_{j+M_i} \hbox{ or vanishes at the origin}.
$$
Of course when $i=n$ the latter possibility obtains. Whenever the case
$m_{j,i}(0)\phi_{j+M_i}$ obtains, we may replace this function by $\phi_{j+M_i}$
and replace $m_{j,n}$ by $m(j,i)(0)m_{j,n}$. Thus we may assume
$$
m_{j,i} \hbox{ is equal to } \phi_{j+M_i} \hbox{ or vanishes at the origin}.
$$


We set $j_n=j$. 
If $m_{j,i}$ vanishes at the origin, we set $j_i=j_n+M_i$. If $m_{j,i} = \phi_{j+M_i}$, we
consider the lacunary decomposition
\[\phi_{j + M_i} = \sum_{j_i \leq j_n+M_i} \phi_{j_i} - \phi_{j_i-1}\ \ \ .\]

The claim then follows by an application of  Theorem \ref{paraproduct}.

\endprf

\end{document}